\documentclass[12pt]{article}
\usepackage{latexsym}

\newtheorem{prop}{Proposition}

\newtheorem{rem}{Remark}
\newtheorem{lem}{Lemma}

\newcommand\bC{\mathbf{C}}
\newcommand\bP{\mathbf{P}}

\newcommand\cC{\mathcal{C}}
\newcommand\cK{\mathcal{K}}
\newcommand\cO{\mathcal{O}}
\newcommand\cV{\mathcal{V}}

\newcommand\qed{{\hfill $\Box$}}

\begin{document}
\begin{center}
Versal deformations and superpotentials for rational curves in smooth
threefolds\\ Sheldon Katz\footnote{email: katz@math.okstate.edu}\\
Department of Mathematics\\ Oklahoma State University\\ Stillwater, OK
74078\\ {\em Dedicated to Herb Clemens on the occasion of his
$60^{\scriptstyle{\mathrm{th}}}$ birthday\/}.
\end{center}

\section{Introduction.}
My interests in the study of curves on Calabi-Yau threefolds have
largely been shaped in two ways: by the tutelage and guidance of Herb
Clemens while I was an instructor at the University of Utah in the
early 80's, and by the interaction between geometry and string theory.
This note is the result of both of these influences.

Herb impressed upon me the importance of the consideration of
elementary examples and development of computational techniques for
these examples when a general theory was not yet apparent.  For
example, in his own work he developed simple computational techniques
for studying the moduli of immersed rational curves 10 years before
Kontsevich's notion of a stable map.  This led to breakthroughs in the
study of the Abel-Jacobi mapping such as the infinite generation of
homological equivalence modulo algebraic equivalence \cite{homalg},
led to his formulation of the famous Clemens conjecture on the
finiteness of rational curves on the generic quintic threefold
\cite{aj} and my proof in degree $d\le7$ \cite{finite} which was only
possible thanks to Herb's generous sharing of his ideas.


Given a curve $C$ in any complex manifold $M$, there is a versal deformation
space $\cK$ parametrizing the local moduli of deformations of $C$ in $M$
\cite{namba}.  Letting $N$ denote the normal bundle of $C$ in $M$, the
result is that there is a neighborhood $U\subset H^0(C,N)$ containing
the origin and a holomorphic obstruction map
\[
F:U\to H^1(C,N)
\]
with $F(0)=0$ and $F'(0)=0$ such that $\cK$ is defined as an analytic space by
\[
\cK=F^{-1}(0).
\]
Note that $\cK$ is locally completely determined by $N$ if either $H^0(C,N)=0$
(in which case $\cK$ is a point) or $H^1(C,N)=0$ (in which case $\cK$ is smooth
of dimension $h^0(C,N)$).

I will focus here on the case when $X$ is a complex threefold and
$C\simeq\bP^1$.  By the previous paragraph, I may as well assume that $C$
contains both deformations and obstruction spaces.  In terms of the
normal bundle $N$ of $C$ in $X$, this means that
$N\simeq\cO(m)\oplus\cO(-n)$, with $m\ge0$ and $n\ge2$.

In the generality considered in this note, the
obstruction map can be computed explicitly.

If $X$ is a Calabi-Yau threefold, then the versal deformation space
has another description arising from string theory, since this space
parametrizes D-branes wrapping holomorphic curves in
$X$.\footnote{Technically, one has to consider curves together with a
flat line bundle; but if the curve is rational, then there is no
distinction.}
The result from physics is that the versal deformation space is a
local description of the space of supersymmetric ground states of a
4-dimensional gauge theory with $N=1$ supersymmetry.  Furthermore, the
gauge theory contains $h^0(C,N)$ massless chiral fields \cite{bsv}.  See
\cite{TASI} for an introduction to the physics of D-branes.

In mathematical terms, the assertion from physics implies the
following description of the versal deformation space: there is a
compact Lie group $G$ acting on $H^0(C,N)$, together with a
$G$-invariant holomorphic function $W$ 
on a neighborhood $U\subset
H^0(C,N)$ of the origin such that
\[
\cK\simeq\mathrm{Crit}(W)/G.
\]
In physics, $G$ arises as the gauge group, and $W$ arises as the 
{\em superpotential\/}.  See \cite{freed} for a description of these
gauge theories written for mathematicians. 

Comparing with the mathematical description of the versal deformation space,
it is clear that these results can only be compatible if $G$ acts trivially.

The main result of this note is a proof in the generality considered
here that
the versal deformation space is defined by the vanishing of the gradient of a
single holomorphic function $W$ on a finite dimensional space.

The superpotential $W$ deserves to be better understood in a precise
mathematical sense, as it connects several areas of geometry.  In
\cite{witten}, $W$ is computed as a function on the (infinite
dimensional) space of $C^\infty$ curves as an extension of the
Abel-Jacobi mapping, and is shown to be holomorphic in a certain sense
(see also \cite{dt}).  In \cite{kklm}, a finite-dimensional expression
for $W$ is given in terms of holomorphic Chern-Simons theory, and its
holomorphicity is illustrated in examples.  It is an interesting
problem under investigation to turn this last description into a
mathematical proof that the versal deformation space is given as the
critical variety of a single holomorphic function on a finite
dimensional space, as is proven here in a special case.

The existence of a superpotential simplifies many calculations in
deformation theory.  For instance, the family of lines on quintic
threefolds has received a lot of attention, and in particular, the
famous 2875 lines can be accounted for by deformation away from the
Fermat quintic \cite{ak}.  It is shown in \cite{bdlr} that much of
this deformation theory can be described simply by adding a
perturbation term to to the superpotential $W=\rho^3\psi^3$ for the
Fermat quintic.  The counting of the 2875 lines as a virtual number on
the Fermat quintic has been described in \cite{ck}.

\bigskip\noindent {\bf Acknowledgements.}  It is a pleasure to thank
H.\ Clemens, B.\ Crauder, M.\ Douglas, A.\ Elezi, R.\ Hartshorne, S.\
Kachru, and C.\ Vafa for helpful conversations.

\section{Versal deformation spaces}
Let $C\simeq\bP^1$ be a smooth rational curve contained in a smooth complex
analytic threefold $X$, with normal bundle $N=\cO(m)\oplus\cO(-n)$.  

In this note, explicit computations will be performed for smooth rational
curves $C$ in threefolds $X$ for which $C$ can be covered by 
coordinate charts $U_0$ and $U_1$ of $X$ such that:
\begin{enumerate}
\item $U_0$ has coordinates $(x,y_1,y_2)$
and $C$ is defined in $U_0$ by $y_1=y_2=0$.
\item $U_1$ has coordinates $(w,z_1,z_2)$ 
and $C$ is defined in $U_1$ by $z_1=z_2=0$.
\item On $C$, the
coordinates $w$ and $x$ are related by $w=x^{-1}$.
\end{enumerate}

\begin{rem}
Examples of such curves abound,  for instance if $X$ is quasi-projective
and $C$ is embedded in $X$ as a line.
\end{rem}

\begin{rem}
Since considerations are purely local, it can be assumed that
$X=U_0\cup U_1$.
\end{rem}

Let $I$ denote the ideal sheaf of $C$ in $X$.  Since $I/I^2\simeq
N^*\simeq \cO(n)\oplus\cO(-m)$, $X$ is determined locally near $C$
modulo $I^2$.  After a possible coordinate change, the gluing map
between the patches therefore takes the form
\begin{equation}
\label{glue}
\begin{array}{ccl}
z_1 &=& x^n y_1 + f(x,y_1,y_2)\\
z_2 &=& x^{-m} y_2 + g(x,y_1,y_2)\\
w &=& x^{-1} + h(x,y_1,y_2),
\end{array}
\end{equation}
where $f,g,h$ are sections of $I^2$, holomorphic in $U_0\cap U_1$
(where in (\ref{glue}), these sections are expressed 
in the coordinates of $U_0$).

\medskip
The goal of this section is to explicitly
construct the versal deformation space of 
these curves, following \cite{namba}.  While the main calculation is
really just an exercise in adapting the more general situation in
\cite{namba}, a fair amount of detail will be given here, both so that
this note will be more or less self-contained, and also because I think
that this construction deserves to be better known.

\medskip
The normal bundle $N$ can be constructed in the usual way by gluing
together trivializations over the $V_i:=U_i\cap C$.  The gluing map is given 
by the matrix
\begin{equation}
\label{F}
F=\left(
\begin{array}{cc}
x^n & 0\\
0 & x^{-m}
\end{array}
\right)
\end{equation}
expressed in the coordinate $x$ on $V_0$.  A section of $N$
on an analytic
open set $V\subset V_i$ will be expressed as a vector of holomorphic
functions $(\phi_1^i,\phi_2^i)$ on $V$.  In particular, if
$\cV=\{V_0,V_1\}$, then the group of \v{C}ech cochains
$C^0(\cV,N)$ can be identified with the space 
\[
C^0(\cV,N)=\left\{(\phi_1^0(x),\phi_2^0(x)),
(\phi_1^1(w),\phi_2^1(w))\right\},
\]
of pairs of vector valued-holomorphic functions in $V_0$ (resp.\
$V_1$).  The convention used here is that the two sections of $N$ defined
above agree on $V_0\cap V_1$ if
\begin{equation}
\label{eqsec}
\left(
\phi_1^1(w), \phi_2^1(w)
\right) = 
\left(
\phi_1^0(x), \phi_2^0(x)
\right)F,
\end{equation}
where $F$ is the transition matrix~(\ref{F}).
It follows that the vector space $H^0(N)\subset C^0(\cV,N)$ is given as the
space of cocycles
\begin{equation}
\label{h0}
\left\{\left(\left(
0, \sum_{i=0}^m a_ix^i\right), 
\left(
0, 
\sum_{i=0}^m a_iw^{m-i}
\right)\right)\right\}.
\end{equation}
The $a_i$ serve as coordinates on $H^0(N)$, and the versal deformation space will
be given explictly by equations in the $a_i$.

In addition, $H^1(N)$ can be given by \v{C}ech cohomology with respect to the
cover $\cV$.  There is the usual coboundary map $\delta:C^0(\cV,N)\to
C^1(\cV,N)$.  I fix the convention throughout that in writing elements
of $C^1(\cV,N)$, i.e.\ sections of $N$ over $V_0\cap V_1$, the trivialization
of $N$ over $V_1$ will be used.
With this convention, the coboundaries are just the sections whose series
expansions have the form
\[
\left(\sum_{i\le -n,\ i\ge 0}b_iw^i, \sum_i c_iw^i\right).
\]
The result is that 
\begin{equation}
\label{h1}
H^1(N)\simeq\left\{\left(\sum_{i=1}^{n-1}b_{-i}w^{-i}, 0\right)\right\}.
\end{equation}
There is also an obvious projection $H:C^1(\cV,N)\to H^1(N)$ defined by the
natural truncation of $\sum_i b_iw^i$:
\begin{equation}
\label{H}
H\left(\sum_{i}b_iw^i, \sum_i c_iw^i\right)=
\left(\sum_{i=1}^{n-1}b_{-i}w^{-i}, 0\right).
\end{equation}

Following \cite{namba}, define the map
\[
K:C^0(\cV,N)\to C^1(\cV,N)
\]
by putting

\begin{equation}
\label{K}
\begin{array}{l}
K\left((\phi_1^0(x),\phi_2^0(x)),
(\phi_1^1(w),\phi_2^1(w))\right)=\\
\qquad \left(\phi_1^1(x^{-1}+h(x,\phi_1^0(x),\phi_2^0(x))-\left(
x^n\phi_1^0(x)+f(x,\phi_1^0(x),\phi_2^0(x))\right),\right.\\
\left.\qquad\phi_2^1(x^{-1}+h(x,\phi_1^0(x),\phi_2^0(x))-\left(
x^{-m}\phi_2^0(x)+g(x,\phi_1^0(x),\phi_2^0(x))\right)\right).
\end{array}
\end{equation}
Note that $K$ has been chosen to have the following property: the equations
\[
y_j=\phi_j^0(x),\ z_k=\phi_k^1(w),\ w=x^{-1}+h(x,\phi_1^0(x),\phi_2^0(x))
\]
patch to define a compact complex curve if and only if
\[
K((\phi_1^0(x),\phi_2^0(x)),
(\phi_1^1(w),\phi_2^1(w)))=0.\footnote{For this to make sense, 
a norm is needed on $C^0(\cV,N)$ and the cochain
$(\phi_1^0(x),\phi_2^0(x)), (\phi_1^1(w),\phi_2^1(w)))$ needs to be
kept sufficiently
small to ensure that this curve actually makes sense in $X$.  As there
is little risk of confusion, I will continue this imprecise practice
throughout and defer to \cite{namba} for the necessary estimates.  In
this way, calculations can be performed using series expansions without
having to worry about convergence issues.}
\]

Continuing to adapt to the notation of \cite{namba}, choose the projection
\[
B:C^1(\cV,N)\to \delta C^0(\cV,N)
\]
defined by
\begin{equation}
\label{B}
B\left(\sum_ib_iw^i,\sum_ic_iw^i\right)=
\left(\sum_{i\le -n,\ i\ge0}b_iw^i,\sum_ic_iw^i\right)
\end{equation}
and the mapping
\[
E_0:\delta C^0(\cV,N)\to C^0(\cV,N)
\]
given by
\begin{equation}
\label{E0}
\begin{array}{l}
E_0\left(\sum_{i\le -n, i\ge0}b_iw^i,\sum_ic_iw^i\right)=\\
\hskip1in
\left((-\sum_{i\le -n}b_ix^{-n-i},-\sum_{i<0}c_ix^{m-i}),
(\sum_{i\ge0} b_iw^i,\sum_{i\ge0}c_iw^i)\right),
\end{array}
\end{equation}
which satisfies $\delta E_0=1$.

Finally define $L:C^0(\cV,N)\to C^0(\cV,N)$ by
\begin{equation}
\label{L}
L(\phi)=\phi+E_0BK\phi-E_0\delta\phi.
\end{equation}
Put 
\begin{equation}
\label{M}
M=\{\phi\mid K(\phi)=0\}.
\end{equation}
As described above, $M$ parametrizes
deformations of $C$, but is way too big to be a versal deformation space.  In
particular, it needs to be cut down so that its tangent space is $H^0(N)$.

\begin{lem} {\rm (Namba)}
$L$ is invertible in a neighborhood of 0, and $L(M)\subset H^0(N)$.
\end{lem}

\bigskip\noindent
{\em Proof:\/} For the first statement, compute $L'(0)=1+E_0B\delta-E_0
\delta=1$.  For the second statement, compute that if $\phi\in M$, then
$\delta L(\phi)=\delta\phi-\delta E_0\delta\phi=\delta\phi-\delta\phi=0$.\qed

\bigskip
Now choose a sufficient small neighborhood $U$ of $0\in H^0(N)$, and
write elements of $U$ as in~(\ref{h0})
\begin{equation}
\label{sa}
s(a)=\left(\left(
0, \sum_{i=0}^m a_ix^i\right), 
\left(
0, 
\sum_{i=0}^m a_iw^{m-i}
\right)\right)
\end{equation}
with $a=(a_0,\ldots,a_m)$ constrained to a suitable neighborhood of 0.  Then 
write
\[
HKL^{-1}\left((0,\sum_{i=0}^m a_ix^i), (0,\sum_{i=0}^m a_i w^{m-i})\right)=
\left(\sum_{i=1}^{n-1}k_i(a_0,\ldots,a_m)w^{-i},0\right),
\]
where $H,K$ and $L$ are respectively given by (\ref{H}), (\ref{K}), and
(\ref{L}).
The main result is
\begin{prop}
The versal deformation space $\cK$ is the analytic space defined in $U$
by the vanishing of the $k_i(a_0,\ldots,a_m),\ 1\le i\le n-1$.
\end{prop}

\bigskip\noindent
{\em Proof:\/} \cite{namba}.\qed

\bigskip
The construction of $\cK$ simplifies somewhat and can be made even more
explicit if $f=f(x,y_2)$ is independent of $y_1$ and is holomorphic in
all of $U_0$, while $g=h=0$.  This is in fact the situation occuring
in the examples in Section~2 of \cite{laufer}.  Such curves will be
called {\em Laufer curves\/}.

In this situation, $L$ is more explicitly given by
\[
L(\phi)=\phi+E_0B\left(\left(-f(x,\phi_2^0),0\right)\right)
\]
which is easy to invert as follows.

\bigskip
Consider the function
\begin{equation}
\label{g}
g(x,a_0,\ldots,a_m)=f(x,\sum_{i=0}^m a_ix^i).
\end{equation}
It has a Taylor expansion
\[
g=\sum_{i=0}^\infty f_i(a_0,\ldots,a_m)x^i,
\]
where the $f_i$ are holomorphic in the $a_i$.  Put
\[
h(x,a_0,\ldots,a_m)=\sum_{i=n}^\infty f_i(a_0,\ldots,a_m)x^{i-n}.
\]
Then compute
\[
L\left((-h,\sum_ia_ix^i),(k_0(a_0,\ldots,a_m),\sum a_iw^{m-i})\right)=
s(a),
\]
where $s(a)$ is the section (\ref{sa}) of $N$.
This means that
\[
L^{-1}\left(s(a)\right)=
\left((-h,\sum_ia_ix^i),(f_0(a_0,\ldots,a_m),\sum a_iw^{m-i})\right).
\]
Direct computation yields
\[
HKL^{-1}\left(\left(
0, \sum_{i=0}^m a_ix^i\right), 
\left(
0, 
\sum_{i=0}^m a_iw^{m-i}
\right)\right)=\left(-\sum_{i=1}^{n-1}f_i(a_0,\ldots,a_m)x^i,0\right).
\]

As in the general case, 
the versal deformation space $\cK$ of $C$ is defined by the equations
\begin{equation}
\label{vdl}
k_1(a_0,\ldots,a_m)=\ldots=k_{n-1}(a_0,\ldots,a_m)=0,\qquad k_i=-f_i.
\end{equation}

\bigskip
In this situation, it is easy to construct a universal family of
curves $\cC$ over $\cK$ explicitly.  The family $\cC\subset \cK\times
X$ is given locally by two equations.  In the patch, $\cK\times U_0$,
$\cC$ is defined by
\[
y_1=-\sum_{j=n}^\infty k_j(a_0,\ldots,a_j)x^{j-n},
\qquad y_2=\sum_{i=0}^m a_i x^i,
\]
and in the patch $\cK\times U_1$, $\cC$ is defined by
\[
z_1=k_0(a_0,\ldots,a_m),\qquad z_2=\sum_{i=0}^m a_iw^{m-i}.
\]
To see that we obtain a family $\cC$,
we only have to check that these glue via (\ref{glue}).  The equation
for $z_2$ is obvious.  In the right hand side of the equation for $z$, we get
$\sum_{i=0}^{n-1}k_i(a_0,\ldots,a_m)x^i$, which simplifies to 
$k_0(a_0,\ldots,a_m)$ when the equations defining $\cK$ are used.

\bigskip\noindent
{\bf Example~1.}  This example is from \cite{laufer}.
\[
\begin{array}{ccl}
z_1 &=& x^3 y_1 + y_2^2+x^2y_2^{2n+1}\\
z_2 &=& x^{-1} y_2
\end{array}
\]
This is the famous example of a contractible curve with normal bundle
$\cO(1)\oplus\cO(-3)$.  Compute from (\ref{g})
\[
g(x,a_0,a_1)=
(a_0+a_1x)^2+x^2(a_0+a_1x)^{2n+1}.
\]
Equation (\ref{vdl}) says to extract the coefficients of $x$ and
$x^2$ and include a minus sign, giving
\[
k_1(a_0,a_1)=-2a_0a_1,\qquad k_2(a_0,a_1)=-(a_1^2+a_0^{2n+1}),
\]
so the versal deformation space is defined by $k_1=k_2=0$.  Note that
the versal deformation space is therefore concentrated at the origin.  In
fact, $C$ is contractible \cite{laufer} and therefore does not deform.

\bigskip
Before stating and proving the main result, a trivial lemma is useful.

\begin{lem}
Consider any analytic function of $x$ and the
$a_i$ of the form
\[
h(x,a_0,\ldots,a_m)=r\left(x,\sum_{i=0}^ma_ix^i\right),
\]
where $r=r(x,y)$ is analytic in $(x,y)$.  Write $h=\sum_{i=0}^\infty
h_i(a_0,\ldots,a_m)x^i$.  Then for all $i,j,k$ such that $0\le j\le m$
and $0\le k+j-i\le m$,
\[
\frac{\partial h_i}{\partial a_j}=\frac{\partial h_k}{\partial a_{k+j-i}}
\]
\end{lem}

\bigskip\noindent
{\em Proof:\/} 
Compute
\begin{equation}
\label{partial}
\frac{\partial h}{\partial a_j} = x^j \frac{\partial g}{\partial y}
\left(x,
\sum_{i=0}^ma_ix^i\right).
\end{equation}
The desired result follows immediately by comparing terms in the series
expansions of $\partial h/\partial a_j$ and $\partial h/\partial a_{k+j-i}$.
\qed

\bigskip
Now suppose that $X$ has trivial canonical bundle, or equivalently,
that $m-n=-2$.  The versal deformation space is described by $m+1$ variables
$a_0,\ldots,a_m$, and is defined by $n-1$ equations.  This restriction
on the canonical bundle implies that there are as many equations as
unknowns in (\ref{vdl}).

\begin{prop}
Suppose that $C\subset X$ is a Laufer curve in a threefold $X$ with trivial
canonical bundle.  Describe the versal deformation space $\cK$ as an analytic
space inside a neighborhood $U$ of 0 in $H^0(N)$ as in (\ref{vdl}).  Then
there exists a single holomorphic function $W$ on $U$ such that $\cK$ is
defined as an analytic space as the subvariety of critical points of $W$.
\end{prop}

\bigskip\noindent
{\em Proof:\/}
More precisely, $W$ will be chosen to satisfy
\[
\frac{\partial W}{\partial a_i} = k_{n-1-i}.
\]
Since the claim is local, it suffices to check the integrability condition
\[
\frac{\partial{k_{n-1-i}}}{\partial a_j}=
\frac{\partial{k_{n-1-j}}}{\partial a_i}.
\]
But this follows immediate from the form of $k_i$ and the Lemma.\qed

\bigskip\noindent
In the Example~1, note that
\[
(k_2,k_1) =-\mathrm{grad}\left(
a_0a_1^2+\frac{a_0^{2n+2}}{2n+2}\right),
\]
so we can take as the superpotential
\[
W=-\left(
a_0a_1^2+\frac{a_0^{2n+2}}{2n+2}\right).
\]

\end{document}